\documentclass{article}

\usepackage{graphicx}
\usepackage{tikz}
\usepackage{amsmath,mathtools,amsfonts}
\usepackage{siunitx}

\newcommand\R{\mathcal{R}}
\newcommand\D{\text{d}}
\newcommand\I{\mathcal{I}}
\renewcommand\L{\mathcal{L}}
\renewcommand\P{\mathcal{P}}
\renewcommand\S{\mathfrak{S}}

\newtheorem{remark}{Remark}

\title{Mixing strategies combined with shape design to enhance productivity of a raceway pond }

\author{Olivier Bernard$^1$, Liudi LU$^{1,2}$, Julien Salomon$^2$}
\date
{%
{\small\textit{$^1$Universit\'e Nice C\^ote d'Azur, Inria BIOCORE, BP93, 06902 Sophia-Antipolis Cedex, France}}\\%
{\small\textit{$^2$Inria ANGE, 75589 Paris Cedex 12, France and Sorbonne Universit\'e, CNRS, Laboratoire Jacques-Louis Lions, 75005 Paris, France}}\\%
}

\begin{document}

\maketitle

\begin{abstract}                % Abstract of not more than 250 words.
This paper focuses on mixing strategies and designing shape of the bottom topographies to enhance the growth of the microalgae in raceway ponds. A physical-biological coupled model is used to describe the growth of the algae. A simple model of a mixing device such as a paddle wheel is also considered. The complete process model was then included in an optimization problem associated with the maximization of the biomass production. The results show that non-trivial topographies can be coupled with some specific mixing strategies to improve the microalgal productivity.
\end{abstract}

\section{Introduction}

Microalgae are photosynthetic organisms whose potential has been proven in the last decade for several biotechnological applications (e.g.~\cite{Chisti2007}). They can be cultivated industrially for cosmetics, pharmaceuticals, food complements and green energy applications~\cite{Wijffels2010}. These micro-organisms can be massively cultivated in closed (e.g.~\cite{Perner2003}) or open photobioreactors. According to the applications, the light can be artificial (for high added value products) or natural.

In this paper, we focus on the cultivation of the algae in a raceway pond. The water is mixed and set in motion in this circular basin by means of a paddle wheel~\cite{Chiaramonti2013}. Studies have shown that the topographies can have an impact on the growth rate of the algae~\cite{Bernard02994713}, whereas mixing the microalgae guarantees that each cell have regularly access to light and necessary nutrients to growth~\cite{Demory2018}. In this paper, we extend the study of~\cite{Bernard02994713} by investigating the optimal combinations between mixing strategies and bottom topographies to enhance algal productivity. We show that non trivial topographies can be obtained associated with some specific mixing strategies.  

The paper is organized as follows. In Section~\ref{sec:model}, we describe the hydrodynamical, the biological and the mixing models and we define the coupled model. Afterwards, we present the optimization problem together with numerical optimization procedure in two frameworks. Eventually, we show some numerical tests to illustrate our approach and study the influence of the topography, depth of the raceway and mixing strategy in the optimization process.

\section{Raceway modeling}\label{sec:model}

The raceway system can be described by a coupling between the hydrodynamics and the dynamics of the photosystems in the algae. The raceway is set in motion by a paddle wheel mixing the algae and modifying their depth and therefore the light flux that they see.

\subsection{Shallow water equations}\label{sec:model_phy}

We model the hydrodynamics of our system by the shallow water equations. More precisely, we consider the smooth steady state solutions of the shallow water equations in a laminar regime, which are governed by the following partial differential equations:
\begin{align}
&\partial_x(h u) = 0, \label{eq:sv1} \\ 
&\partial_x(h u^2+g\frac{h^2}2) = -g h \partial_x z_b, \label{eq:sv2}
\end{align}
where $h$ stands for the water elevation, $u$ represents the horizontal averaged velocity of the water, the constant $g$ is the gravitational acceleration, and $z_b$ defines the topography. The free surface $\eta$ is defined by $\eta := h + z_b$ and the averaged discharge $Q=h u$. A schematic representation of this system is given in Fig.~\ref{fig:P}.
\usetikzlibrary{decorations.pathreplacing}
\begin{figure}[hptb]
\begin{center}
\begin{tikzpicture}
\node at (0,-0.2)[anchor=north] {0};
\node at (3,-0.2)[anchor=north] {$L$};
\node at (4,-0.2)[anchor=north] {0};
\node at (7,-0.2)[anchor=north] {$L$};
\draw [fill,gray](0,-0.1)--(0,0.2) sin (0.75,0) cos (1.5,0.2) sin (2.25,0.4) cos (3,0)--(3,-0.1);
\draw [fill,gray](4,-0.1)--(4,0.2) sin (4.75,0) cos (5.5,0.2) sin (6.25,0.4) cos (7,0)--(7,-0.1);
\draw [dotted,thick](0,1.2) sin (0.75,1) cos (1.5,1.2) sin (2.25,1.4) cos (3,1);
\draw [dotted,thick](4,1.2) sin (4.75,1) cos (5.5,1.2) sin (6.25,1.4) cos (7,1);
\draw [dotted,thick](0,2.2) sin (0.75,2) cos (1.5,2.2) sin (2.25,2.4) cos (3,2);
\draw [dotted,thick](4,2.2) sin (4.75,2) cos (5.5,2.2) sin (6.25,2.4) cos (7,2);
\draw [thick](0,3.2) sin (0.75,3) cos (1.5,3.2) sin (2.25,3.4) cos (3,3);
\draw [thick](4,3.2) sin (4.75,3) cos (5.5,3.2) sin (6.25,3.4) cos (7,3);
\draw [dotted,->](3,0) -- (3.9,3.2);
\draw [dotted,->](3,1) -- (3.9,0.2);
\draw [dotted,->](3,2) -- (3.9,1.2);
\draw [dotted,->](3,3) -- (3.9,2.2);
\node at (3.5,4) (P){$P$};
\draw [->,thick](P)-- (3.5,3);
\draw [->,thick](0,-0.1) -- (0,4.5);
\draw [thick](0,-0.1) -- (7,-0.1);
\draw [dashed,->](0,3.2) -- (7.5,3.2);
\node at (-0.4,4.5) {$z$};
\node at (7.5,3.4) {$x$};
\node at (-0.2,3.2) {0};
\node at (7.5,2.8) {$z_1=z_{\sigma(4)}$};
\node at (7.5,2.2) {$z_2=z_{\sigma(1)}$};
\node at (7.5,0.8) {$z_3=z_{\sigma(2)}$};
\node at (7.5,0.2) {$z_4=z_{\sigma(3)}$};
\draw [->](0.8,1.5) -- (1,1.5);
\draw [->](0.8,2.5) -- (1,2.5);
\draw [->](0.8,0.5) -- (1,0.5);
\draw [->](6,1.8) -- (6.2,1.8);
\draw [->](6,2.8) -- (6.2,2.8);
\draw [->](6,0.8) -- (6.2,0.8);
\draw [dashed,->](0.75,3.2)--(0.75,3);
\draw [dashed,->](2.25,3.2)--(2.25,0.4);
\draw [dashed,<-](4.75,3)--(4.75,0);
\node at (0.75,3.4) {$\eta(x)$}; 
\node at (2.7,1.7) {$z_b(x)$}; 
\node at (4.4,1.5) {$h(x)$}; 
\node at (1.4,1.7) {$u(x)$}; 
\node at (5.6,1.7) {$u(x)$}; 
\draw [thick,->](2.25,4)--(2.25,3.4);
\draw [thick,->](6.25,4)--(6.25,3.4);
\draw [thick,->](1.75,4)--(1.75,3.3);
\draw [thick,->](5.75,4)--(5.75,3.3);
\draw [thick,->](2.75,4)--(2.75,3.2);
\draw [thick,->](6.75,4)--(6.75,3.2);
\node at (2.25,4.25) {$I_s$};
\node at (6.25,4.25) {$I_s$};
\end{tikzpicture}
\end{center}
\caption{Representation of the hydrodynamic model with an example of mixing device ($P$). Here, $P$ corresponds to the cyclic permutation $\sigma= (1 \ 2 \ 3\ 4)$.}
\label{fig:P}
\end{figure}
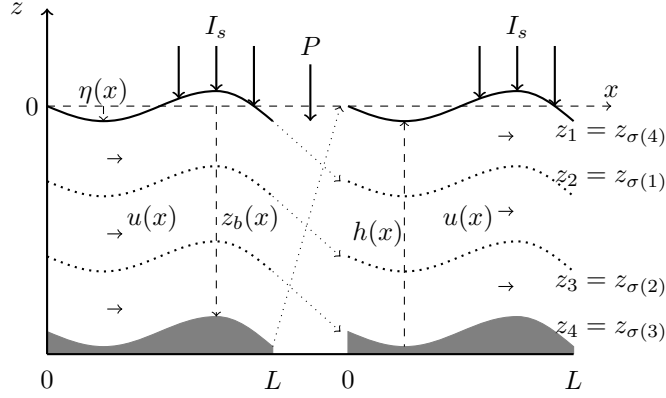

The $z$ axis represents the vertical direction and the $x$ axis represents the horizontal direction. Besides, $I_s$ represents the light intensity at the free surface (assumed to be constant).

The Froude number for the steady state is defined by $Fr=u/\sqrt{g h}$. The situation $Fr<1$ corresponds to the subcritical case (i.e. the flow regime is fluvial) while $Fr>1$ is to the supercritical case (i.e. the flow regime is torrential). In the following studies, we limit ourselves to the subcritical case. Following the procedure from~\cite{Bernard02994713}, the shallow water equations~\eqref{eq:sv1}-\eqref{eq:sv2} can be transformed to: 
\begin{align}
u &= \frac{Q_0}h,\label{eq:u}\\
z_b &= \frac{M_0}g - \frac{Q_0^2}{2g h^2} - h, \label{eq:zb}
\end{align}
where $Q_0$ and $M_0$ are two positive constants.

\subsection{Lagrangian trajectories and light intensity} 

Incompressibility and a kinematic condition at bottom can then be used to obtain the following equation of the Lagrangian trajectory of the algae~\cite[Equation (12)-(15)]{Bernard02994713}.
\begin{equation}\label{eq:zx}
z(x) = \eta(x) + \frac{u(0)}{u(x)}(z(0) - \eta(0)),
\end{equation}
where $z(0)$ is the initial position of the algae. 

To obtain the light intensity $I$ at depth $z$, we use the Beer-Lambert law to describe the light attenuation:
\begin{equation}\label{eq:Beer}
I(x,z) = I_s\exp\big(-\varepsilon(\eta(x) - z)\big),
\end{equation}
where $\varepsilon$ is the extinction coefficient. Using~\eqref{eq:zx} in~\eqref{eq:Beer}, we find the light intensity captured by the algae following the trajectories $z(x)$ in one lap of the raceway
\begin{equation*}
I(x,z) = I_s\exp\big(-\varepsilon \frac{u(0)}{u(x)}(\eta(0) - z(0))\big).
\end{equation*}
Given initial conditions $h(0),z_b(0)$, then $\eta(0)=h(0)+z_b(0)$ and $u(0)$ from~\eqref{eq:u}, we see that computing $I$ on a trajectory only requires to know the initial position $z(0)$ and the velocity $u(x)$, i.e., according to~\eqref{eq:u}, $h(x)$ which is obtained by solving~\eqref{eq:zb}.

\subsection{Han model}

The dynamics of the light harvesting complexes in the chloroplasts is controlled by the amount of light perceived by the algal cells. They can be described by the Han model~\cite{Han2001}, in which the reaction centers are assumed to have three different states: open and ready to harvest a photon ($A$), closed while processing the absorbed photon energy ($B$), or inhibited if several photons have been absorbed simultaneously leading to an excess of energy ($C$). The evolution of the state $A, B, C$ satisfies the following dynamical system
\begin{equation*}
\left\{
\begin{array}{lr}
\dot{A} = -\sigma I A + \frac B{\tau},\\
\dot{B} =  \sigma I A - \frac B{\tau} + k_r C - k_d\sigma I B,\\
\dot{C} = -k_r C + k_d \sigma I B,
\end{array}
\right.
\end{equation*}
where $A, B, C$ are the relative frequencies of the three possible states which satisfy 
\begin{equation*}
A+B+C=1,
\end{equation*}
the coefficients $\sigma$, $\tau$, $k_r$ and $k_d$ represent the specific photon absorption, the turnover rate, the photosystem repair rate and the damage rate, respectively. As shown in~\cite{Lamare2018}, one can use a fast-slow approximation and singular perturbation theory to reduce this system to a single evolution equation:
\begin{equation}\label{eq:evolC}
\dot{C} = -\alpha(I) C + \beta(I),
\end{equation}
where
\begin{equation*}
\alpha(I) = k_d\tau \frac{(\sigma I)^2}{\tau \sigma I+1} + k_r \text{ and } \beta(I) =  \alpha(I) - k_r.
\end{equation*}
Then following~\cite{Bernard02994713}, we obtain a time-free reformulation of~\eqref{eq:evolC}, namely
\begin{equation}\label{eq:Cx}
C' = \frac{-\alpha(I) C + \beta(I)}{u},
\end{equation}
where all the functions on the right-hand side only depend on $x$. The net specific growth rate is then obtained by balancing photosynthesis and respiration, which gives
\begin{equation}\label{eq:mu}
\mu(C,I) = \frac{-\gamma(I)C + \zeta(I)}{u},
\end{equation}
where
\begin{equation*}
\zeta(I) = \frac{k\sigma I}{\tau \sigma I+1}-R \text{ and } \gamma(I) = \zeta(I) + R.
\end{equation*}
Here $k$ stands for a factor that links received energy and growth rate. The term $R$ represents the respiration rate. The average net specific growth rate over the domain is then defined from~\eqref{eq:mu} by
\begin{equation}\label{eq:mucontinuous}
\bar{\mu} := \frac 1L\int_0^L\frac 1{h(x)}\int_{z_b}^{\eta} \mu\big(C(x,z), I(x,z)\big) \D z \D x.
\end{equation}
This will be the principle function of our following studies. 

\begin{remark}\label{rm:nonhomo}
The dynamic of the biomass concentration is derived from~\eqref{eq:mu}:
\begin{equation*}
\dot{X} = \bar{\mu}(C,I)X - D X,
\end{equation*}
where $D$ is the dilution rate. The extinction coefficient $\varepsilon$ depends on $X$ as follows. The system is perfectly mixed, then the concentration is homogeneous so $\varepsilon$ is a constant independent of space. In general, the extinction coefficient $\varepsilon$ is an affine function of the biomass $X$ (see~\cite{Martinez201811})
\begin{equation}\label{eq:eps}
\varepsilon(X) = \alpha_0 X + \alpha_1,
\end{equation}
where $\alpha_0>0$ stands for the specific light extinction coefficient of the microalgae specie and $\alpha_1$ defines the background turbidity that summarizes the light absorption and diffusion due to all non-microalgae components. 

We assume that the algal biomass $X$ in the raceway is controlled at a concentration which meets the so-called compensation condition (see~\cite{Masci2010,Grognard2014}). This condition means that photosynthesis equilibrates respiration in the darkest layer of the raceway. In other terms, at steady-state, the growth rate $\mu$ at the (average) bottom depth $\bar z_b$ is 0, i.e.,
\begin{equation*}
-\gamma(I_{\bar z_b})\frac{\beta(I_{\bar z_b})}{\alpha(I_{\bar z_b})} + \zeta(I_{\bar z_b}) = 0.
\end{equation*}
Solving the above equation provides a value of $I_{\bar z_b}$ thus of the extinction coefficient $\varepsilon(X)$, and finally of the biomass concentration $X$. In the sequel we assume that an appropriate control strategy maintains the biomass around this value by playing on the dilution rate $D$.
\end{remark}

\subsection{Vertical discretization of the system}

In order to compute numerically~\eqref{eq:mucontinuous}, let us consider a uniform vertical discretization of the initial position $z(0)$ for $N_z$ cells:
\begin{equation*}
z_n(0) = \eta(0) - \frac {n-\frac12}{N_z}h(0), \quad n=1,\ldots,N_z. 
\end{equation*}
From \eqref{eq:zx}, we obtain
\begin{equation*}
z_{n}(x) - z_{n+1}(x) = \frac 1{N_z}h(x), \quad n=1,\ldots,N_z, 
\end{equation*}
meaning that the cell distribution remains uniform along the trajectories. To simplify notations, we write $I_n(x)$ instead of $I(x,z_n)$ hereafter. 

Let $C_n(x)$ (resp. $I_n(x)$) the photo-inhibition state (resp. the light intensity) associated with the trajectories $z_n(x)$. Then the semi-discrete average net specific growth rate of~\eqref{eq:mucontinuous} can be defined by
\begin{equation}\label{eq:mudiscret}
\begin{split}
\bar \mu_{\Delta} :&= \frac 1L \int_0^L \frac 1{h(x)}\frac 1{N_z}h(x)\sum_{n=1}^{N_z}\mu( C_n(x),I_n(x) ) \D x\\
&=\frac 1L \frac 1{N_z}\sum_{n=1}^{N_z} \int_0^L\mu( C_n(x),I_n(x) ) \D x.
\end{split}
\end{equation}

%% There are a number of predefined theorem-like environments in
%% ifacconf.cls:
%%
%% \begin{thm} ... \end{thm}            % Theorem
%% \begin{lem} ... \end{lem}            % Lemma
%% \begin{claim} ... \end{claim}        % Claim
%% \begin{conj} ... \end{conj}          % Conjecture
%% \begin{cor} ... \end{cor}            % Corollary
%% \begin{fact} ... \end{fact}          % Fact
%% \begin{hypo} ... \end{hypo}          % Hypothesis
%% \begin{prop} ... \end{prop}          % Proposition
%% \begin{crit} ... \end{crit}          % Criterion

\subsection{Paddle wheel modeling}

Recent studies have shown that the paddle wheel played a key role in a raceway ponds system~\cite{Chiaramonti2013,Demory2018}, where the paddle wheel set this hydrodynamic-biologic coupling system in motion. At the same time, it modifies the elevation of the algae passing through it, and thus giving successively access to light to all the population. This mixing device has been studied in~\cite{Bernard02970756}, where a flat topography has been considered and the mixing procedure is assumed to be perfect, meaning that at each new lap, the algae at depth $z_n(0)$ are entirely transferred into the position $z_{\sigma(n)}(0)$ when passing through the mixing device. In the current study, we still assume that we can design a mixing setup achieving an ideal rearrangement of trajectories, and we consider the case when the topography is no longer flat.

We denote by $\P$ the set of permutation matrices of size $N_z\times N_z$ and by $\S_{N_z}$ the associated set of permutations of $N_z$ elements. This model is depicted schematically on an example in Fig.~\ref{fig:P}.

\subsection{Periodic regime}

We assume that the state $C$ is $KL$-periodic in the sense that after $K$ times of passing the device ($P$),  $C^K(0)=C(0)$. A natural choice for $K$ is the order of the permutation $P$. 

Following arguments similar to that in~\cite[Proposition 1, Lemma 1]{Bernard02970756}, it can be proved that if the system is $KL$-periodic, then it is $L$-periodic. Hence, the average growth rate $\bar \mu$ for $K$ laps equals to the average growth rate $\bar \mu$ for a single lap. This will help us in simplifying the formulation of the optimization problem considered in the next section. In addition, the computations to solve the optimization problem will be reduced, since the CPU time required to assess the productivity gain of a permutation will not depend on its order.
\begin{remark}
In the setting presented in~\cite{Bernard02970756}, when the system is assumed to be periodic $C(0)$, hence $C$ depends on the permutation matrix $P$. In the current study, the state $C$ will also depend on the permutation matrix $P$ that we denote $C^P$ hereafter.
\end{remark}

\section{Optimization}

In this section, we define the optimization problem associated with our biological-hydrodynamical-permutation model. As mentioned in Section~\ref{sec:model_phy}, a given smooth topography $z_b$ corresponds to a unique water elevation $h$ under the assumption that flow remains in a subcritical regime. On the other hand, since we consider a 1D framework, the volume of our system is simply given by
\begin{equation*}
V = \int_0^L h(x) \D x.
\end{equation*}
Since this quantity plays an important role in raceway design, and need to be easily handled. Therefore, we choose to parameterize the water elevation $h$. Given an optimal parameter $a^*$, the associated optimal topography can be determined by means of~\eqref{eq:zb}. An example of parameterization consists in writing $h$ as a truncated Fourier series
\begin{equation}\label{eq:fourier}
h(x,a) = a_0 + \sum_{m=1}^M a_m \sin (2m\pi \frac x L),
\end{equation}
Such parameterization allows us to control the volume since $V=a_0L$.

For simplicity of notation, we omit $x$ in the notation and rather denote explicitly that the functions depend on $a$.

\subsection{Optimization problem for constant reactor volume}

In this section, we assume that the volume of the reactor is constant. Such situation can be obtained, e.g., by using the parameterization~\eqref{eq:fourier} with a fixed $a_0$. Consider then a vector $a:=[a_1,\cdots,a_M]\in \R^{M}$, which will be the variable to be optimized. The objective function is then defined from~\eqref{eq:mudiscret} by
\begin{equation}\label{eq:muPa}
\bar \mu_{\Delta}^P(a) = \frac 1{LN_z}\sum_{n=1}^{N_z}\int_0^L\frac{-\gamma(I_n(a))C^P_n + \zeta(I_n(a))}{u(a)}   \D x,
\end{equation}
where $C^P_n$ satisfies the following parameterized version of~\eqref{eq:Cx} with a periodic condition
\begin{equation}\label{eq:cond}
\left\{
\begin{array}{lr}
{C^P_n}' + \frac{\alpha(I_n(a))}{u(a)}C^P_n & =  \frac{\beta(I_n(a))}{u(a)}\\
PC^P_n(L) & = C^P_n(0).
\end{array}
\right.
\end{equation}
Our optimization problem then reads:

\textit{Find a permutation matrix $P_{\max}$ and a parameter vector $a^*$ solving the maximization problem:}
\begin{equation*}
\max_{P\in \P} \max_{a\in\R^M}\bar \mu_{\Delta}^P(a).
\end{equation*}

\subsection{Optimization procedure for constant reactor volume}

For a given permutation matrix $P\in \P$, the Lagrangian of~\eqref{eq:muPa} can then be written by
\begin{equation*}
\begin{split}
\L^P(C,p,a) = \frac 1{LN_z}&\sum_{n=1}^{N_z}\int_0^L \frac{-\gamma(I_n(a))C^P_n +\zeta(I_n(a))}{u(a)} \D x\\
-&\sum_{n=1}^{N_z}\int_0^L p^P_n\big({C^P_n}'+ \frac{\alpha(I_n(a))C^P_n-\beta(I_n(a))}{u(a)}\big) \D x
\end{split}
\end{equation*}
where $p^P_n$ is the Lagrange multiplier associated with the constraint~\eqref{eq:cond}.

The optimality system is obtained by cancelling all the partial derivatives of $\L^P$. Differentiating $\L^P$ with respect to $p^P_n$ and equating the resulting expression to zero gives~\eqref{eq:cond}. Integrating the terms $\int p^P_n {C^P_n}'\D x$ on the interval $[0,L]$ by parts enables to differentiate $\L^P$ with respect to $C^P_n$ and $C^P_n(L)$. Equating the result to zeros gives rise to
\begin{equation}\label{eq:adjoint}
\left\{
\begin{array}{lr}
{p_n^P}' - p_n^P \frac{\alpha(I_n(a))}{u(a)} -\frac 1{LN_z}\frac{\gamma(I_n(a))}{u(a)} &= 0 \\
p_n^P(L) - p_n^P(0)P &= 0.\\
\end{array}
\right.
\end{equation}
Given a vector $a$, let us still denote by $C_n^P, p_n^P$ the corresponding solutions of~\eqref{eq:cond} and~\eqref{eq:adjoint}. The gradient $\nabla \bar \mu^P_{\Delta}(a)$ is obtained by
\begin{equation*}
\nabla \bar \mu^P_{\Delta}(a) = \partial_{a}\L^P,
\end{equation*}
where
\begin{equation*}
\begin{split}
\partial_{a}\L^P = \frac 1{LN_z}&\sum_{n=1}^{N_z}\int_0^L \frac{-\gamma'(I_n(a))C^P_n+\zeta'(I_n(a))}{u(a)}\partial_a I_n(a)\D x \\
-\frac 1{LN_z}&\sum_{n=1}^{N_z}\int_0^L \frac{-\gamma (I_n(a))C^P_n + \zeta(I_n(a))}{u^2(a)}\partial_a u(a)\D x \\
+&\sum_{n=1}^{N_z}\int_0^L p_n^P\frac{-\alpha'(I_n(a))C^P_n + \beta'(I_n(a))}{u(a)}\partial_a I_n(a)\D x \\
-&\sum_{n=1}^{N_z}\int_0^L p_n^P\frac{-\alpha(I_n(a))C^P_n + \beta(I_n(a))}{u^2(a)}\partial_a u(a) \D x.
\end{split}
\end{equation*}

\subsection{Optimization problem for variable reactor volume}

In this section, we focus on the case where the reactor volume can also vary. As we have mentioned in Remark~\ref{rm:nonhomo}, we apply an extra assumption to determine $X$ as a function of the volume. Therefore, we apply the parameterization~\eqref{eq:fourier} and follow the computations in~\cite{Bernard02970756} to determine $X$. Such parameterization allows us to control the biomass $X$ and the volume of the system $V$ by using an extra parameter $a_0$. Note that the larger $N$, the less valid is our hydrodynamic model, see Section~\ref{sec:model_phy}, where a smooth topography is assumed to guarantee a laminar regime. Hence, limit situations where $N\rightarrow+\infty$ are not considered in what follows. Since the biomass $X$ and $V$ the volume of the system vary with the parameter $a_0$, maximizing areal productivity is a more relevant target. For a given biomass concentration $X$, the productivity per unit of surface is given by:
\begin{equation}\label{eq:picontinous}
\Pi = \bar \mu X \frac VS,    
\end{equation}
where $S$ presents the ground surface of the raceway pond. From~\cite[Appendix C]{Bernard02994713}, we have
\begin{equation*}
X\frac VS = \alpha_2 - \alpha_3 a_0 \text{ with } \alpha_2=\frac{1}{\alpha_0}\ln(\frac{I_s}{I_{\bar z_b}}) \text{ and } \alpha_3=\frac{\alpha_1}{\alpha_0},
\end{equation*}
where $\alpha_0, \alpha_1$ are given in~\eqref{eq:eps}.

Consider the extend parameter vector $\tilde a:=[a_0, a]\in \R^{M+1}$. From~\eqref{eq:picontinous}, the objective function is given by
\begin{equation}\label{eq:piPa}
\Pi_{\Delta}^P(\tilde a):= \bar \mu_{\Delta}^P(\tilde a) (\alpha_2 - \alpha_3 a_0),
\end{equation}
where $\bar \mu_{\Delta}^P$ is defined in~\eqref{eq:muPa}.

The corresponding optimization problem reads: 

\textit{Find a permutation matrix $P_{\max}$ and a parameter vector $\tilde a^*$ solving the maximization problem:}
\begin{equation*}
\max_{P\in \P}\max_{\tilde a\in \R^{M+1}} \Pi_{\Delta}^P(\tilde a).
\end{equation*}

\subsection{Optimization procedure for variable reactor volume}

Let us denote by $\tilde \L^P$ the Lagrangian associated to~\eqref{eq:piPa}. We follow the same optimization procedure presented in Section 3.2. Note that an extra element appears in this gradient, which is the partial derivative of $\tilde \L^P$ with respect to the variable $a_0$. More precisely, we have $\nabla \Pi_{\Delta}^P(\tilde a) = [\partial_{a_0}\tilde \L^P,\partial_{a}\tilde \L^P]$ where
\begin{equation*}
\begin{split}
\partial_{a_0}\tilde \L^P = \frac {\alpha_2 - \alpha_3 a_0}{LN_z}&\sum_{n=1}^{N_z}\int_0^L \frac{-\gamma'(I_n(\tilde a))C^P_n+\zeta'(I_n(\tilde a))}{u(\tilde a)}\partial_{a_0} I_n(\tilde a)\D x \\
-\frac {\alpha_2 - \alpha_3 a_0}{LN_z}&\sum_{n=1}^{N_z}\int_0^L \frac{-\gamma (I_n(\tilde a))C^P_n + \zeta(I_n(\tilde a))}{u^2(\tilde a)}\partial_{a_0} u(\tilde a)\D x \\
-\frac {\alpha_3}{LN_z}&\sum_{n=1}^{N_z}\int_0^L \frac{-\gamma (I_n(\tilde a))C^P_n + \zeta(I_n(\tilde a))}{u(\tilde a)}\D x \\
+&\sum_{n=1}^{N_z}\int_0^L p_n^P\frac{-\alpha'(I_n(\tilde a))C^P_n + \beta'(I_n(\tilde a))}{u(\tilde a)}\partial_{a_0} I_n(\tilde a)\D x \\
-&\sum_{n=1}^{N_z}\int_0^L p_n^P\frac{-\alpha(I_n(\tilde a))C^P_n + \beta(I_n(\tilde a))}{u^2(\tilde a)}\partial_{a_0} u(\tilde a) \D x,
\end{split}
\end{equation*}
and $\partial_{a}\tilde \L^P$ is very similar to $\partial_{a}\L^P$, by adding an extra product of $(\alpha_2 - \alpha_3 a_0)$ with the first-two integrals in $\partial_{a}\L^P$ and change $a$ to $\tilde a$ in all the parameter-depend functions.

\section{Numerical results}

In this section, we present some numerical results derived from the optimization procedure presented in the previous section. Note that for a given vertical discretization number $N_z$, we need to test all the permutation matrices in the set $\P$, which means $N_z!$ possible cases. Since an optimization problem must be solved for each permutation matrix the problem is highly computational demanding and  we consider a small value of $N_z$ for which the problem is solvable in reasonable time. Regarding the parameterization of $h$, we use truncated Fourier series presented in~\eqref{eq:fourier}.

\subsection{Numerical solver}

We introduce a supplementary space discretization with respect to $x$ to solve numerically our optimization problem. Let us take a space increment $\Delta x$, set $N_x =[L/\Delta x]$ and $x^{n_x}=n_x\Delta x$ for $n_x=0,\ldots,N_x$. We choose to apply the Heun's scheme to compute $C^P_n$ via~\eqref{eq:cond}. Following a first-discretize-then-optimize strategy, we get that the Lagrange multiplier $p^P_n$ is also computed by a Heun's type scheme via~\eqref{eq:adjoint}. We use the \textit{fmincon} solver in MATLAB to solve the optimization problem with the subcritical constraint.

\subsection{Parameter for the models}

The spatial increment is set to $\Delta x=\SI{0.01}{m}$ such that the convergence of the numerical scheme has been ensured, and we take the averaged discharge $Q_0 = \SI{0.04}{m^2\cdot s^{-1}}$, and $z_b(0) = \SI{-0.4}{m}$ to stay in standard ranges for a raceway pond. The free-fall acceleration is set to be $g=\SI{9.81}{m\cdot s^{-2}}$. All the numerical parameters values for Han's model are taken from~\cite{Grenier2020} and given in table~\ref{tab:hanpara}.
\begin{table}[htpb]
\caption{Parameter values for Han Model}
\label{tab:hanpara}
\begin{center}
\begin{tabular}{|c|c|c|}
\hline
$k_r$  & 6.8 $10^{-3}$ & s$^{-1}$\\
\hline
$k_d$ & 2.99 $10^{-4}$  & -\\
\hline
$\tau$ & 0.25 & s\\
\hline
$\sigma$ & 0.047 & m$^2 \cdot(\mu$ mol)$^{-1}$\\
\hline
$k$  & 8.7 $10^{-6}$ & -\\
\hline
$R$ &  1.389 $10^{-7}$ & s$^{-1}$\\
\hline
\end{tabular}
\end{center}
\end{table}
For fixed volume, we assume that only $1\%$ of light can be captured by the cells at the bottom of the raceway, under our parameterization, the light extinction coefficient $\varepsilon$ can be computed by
\begin{equation*}
\varepsilon = (1/a_0)\ln(1/1\%).
\end{equation*}
For varying volume, the specific light extinction coefficient of the microalgae specie $\alpha_0=\SI{0.2}{m^2\cdot g}$C and the background turbidity coefficient $\alpha_1=\SI{10}{m^{-1}}$, these are taken from~\cite{Martinez201811}. Besides, $I_s=\SI{2000}{\mu mol\cdot m^{-2}.s^{-1}}$ which corresponds to the order of magnitude of the maximum light intensity in summer in the south of France.

\subsection{Numerical tests}

We present some results for both constant and non constant volume. We will show the optimal permutation matrix and the associated shape of the topography for these two frameworks and compare the gain with a standard system.

\subsubsection{Constant volume}

\begin{figure}[hptb]
\begin{center}
\includegraphics[scale=0.3]{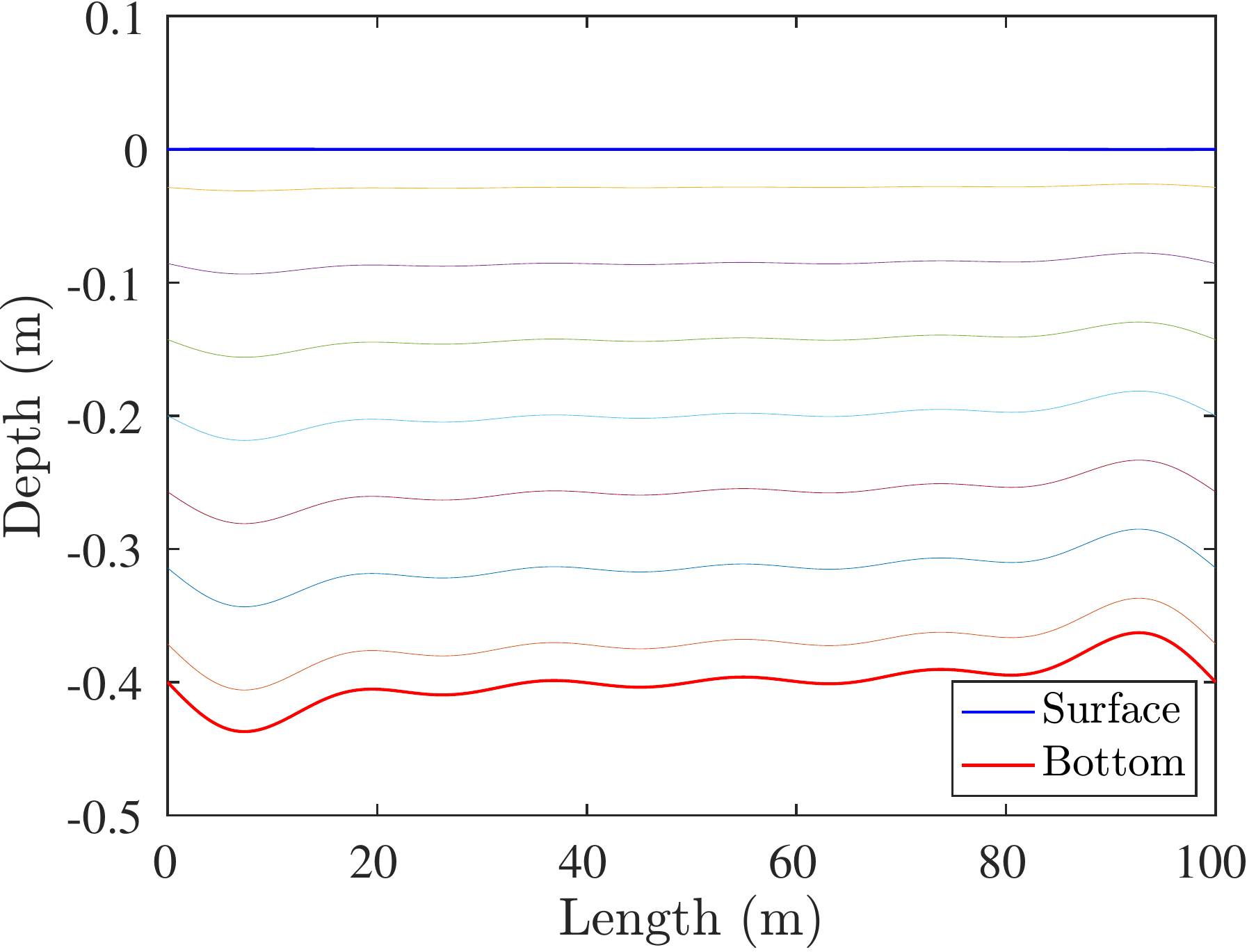}
\end{center}
\caption{The optimal topography and the associated trajectories for the permutation matrix~\eqref{eq:matL100}.}
\label{fig:top1}
\end{figure}

The first test is dedicated to study the optimal permutation matrix and the associated shape of the topography when the volume is constant. To evaluate the efficiency of the corresponding mixing strategy, let us define:
\begin{align}
&r_1:=\frac{\bar \mu_{\Delta}^{P_{\max}}(a^*)-\bar \mu_{\Delta}^{P_{\max}}(0)}{\bar \mu_{\Delta}^{P_{\max}}(0)}, \label{eq:r1}\\
&r_2:=\frac{\bar \mu_{\Delta}^{P_{\max}}(a^*)-\bar \mu_{\Delta}^{\I_{N_z}}(0)}{\bar \mu_{\Delta}^{\I_{N_z}}(0)}.\label{eq:r2}
\end{align}
Here $r_1$ defines the gain of the optimal permutation strategy with the optimal topography compare to the optimal permutation strategy with a flat topography, and $r_2$ defines the gain of the optimal permutation strategy with the optimal topography compare to no permutation strategy with a flat topography. We fix the parameter $a_0(=h(0;a))=\SI{0.4}{m}$ to stay in a standard raceway pond range. Let us take the truncated Fourier number $M = 5$ and the vertical discretization number $N_z = 7$. The initial guess of the vector $a$ is set to be 0, which corresponds to a flat topography. Let us start the test with a given raceway pond length $L=\SI{100}{m}$. The optimal matrix $P_{\max}$ is given in~\eqref{eq:matL100} and the associated optimal topography is presented in Fig.~\ref{fig:top1}. A non flat topography associated with a non trivial permutation matrix has been observed. The optimal value of $a$ is $[0.0123, 0.0119, 0.0097, 0.0080, 0.0067]$. In particular, this optimal matrix corresponds to the optimal matrix obtained with a flat topography under the same parameter settings~\cite[Equation 15]{Bernard02970756}. Besides, in this case $r_1=0.148\%$ and $r_2=1.070\%$.
\begin{equation}\label{eq:matL100}
P_{\max}=\begin{pmatrix}
0 & 1 & 0 & 0 & 0 & 0 & 0  \\
0 & 0 & 0 & 1 & 0 & 0 & 0  \\
0 & 0 & 0 & 0 & 0 & 1 & 0  \\
0 & 0 & 0 & 0 & 0 & 0 & 1  \\
0 & 0 & 0 & 0 & 1 & 0 & 0  \\
0 & 0 & 1 & 0 & 0 & 0 & 0  \\
1 & 0 & 0 & 0 & 0 & 0 & 0 
\end{pmatrix}
\end{equation}

\begin{figure}[hptb]
\begin{center}
\includegraphics[scale=0.3]{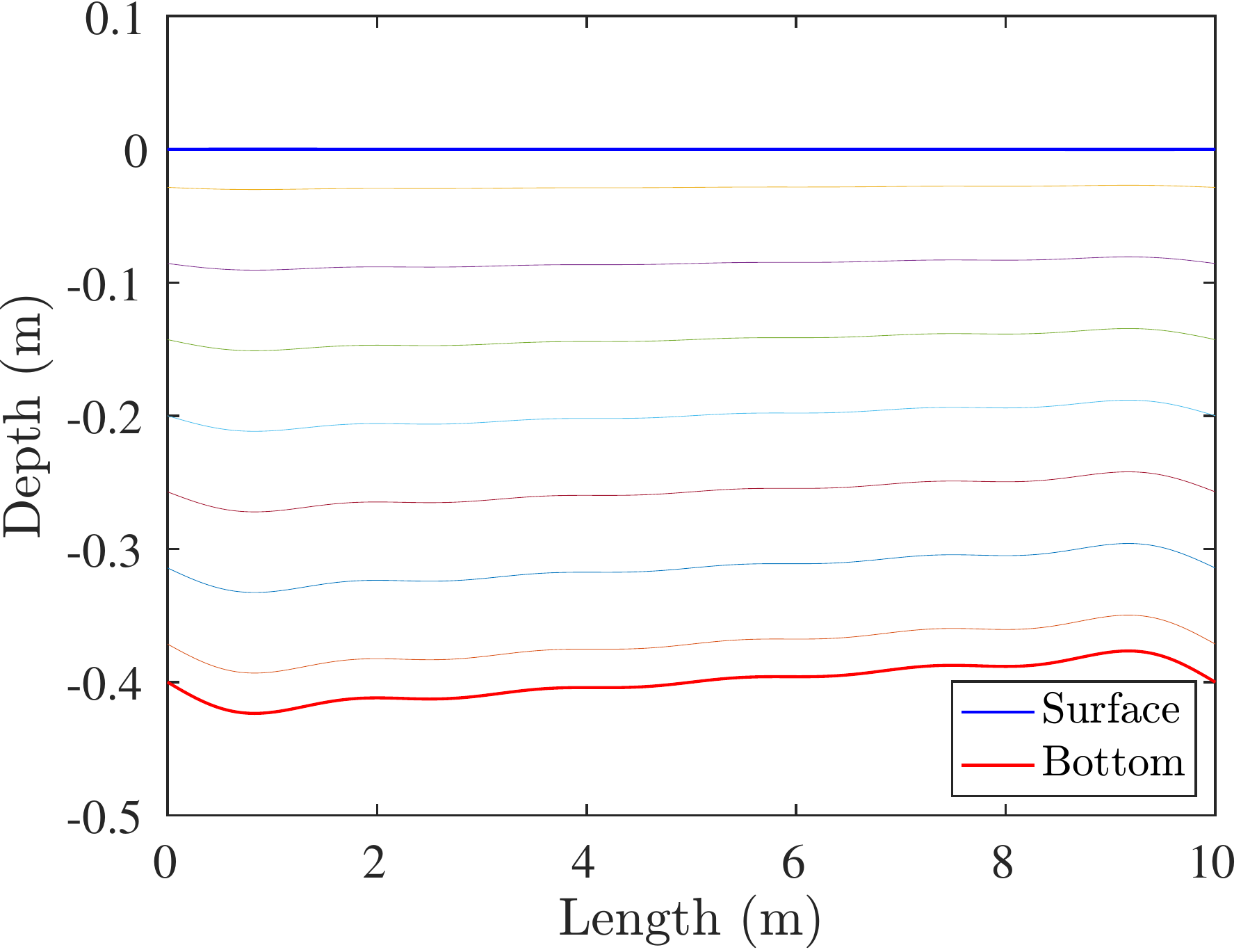}
\end{center}
\caption{The optimal topography and the associated trajectories for the permutation matrix~\eqref{eq:matL10}.}
\label{fig:top2}
\end{figure}

We then proceed a test for the length $L=\SI{10}{m}$. The optimal matrix $P_{\max}$ is given in~\eqref{eq:matL10} and the associated optimal topography is presented in Fig.~\ref{fig:top2}. The optimal value of $a$ is $[0.0147,0.0074,0.0050,0.0037,0.0030]$, the two ratios defined in~\eqref{eq:r1}-\eqref{eq:r2} are $r_1=0.089\%$ and $r_2=3.542\%$.
\begin{equation}\label{eq:matL10}
P_{\max}=\begin{pmatrix}
1 & 0 & 0 & 0 & 0 & 0 & 0  \\
0 & 0 & 1 & 0 & 0 & 0 & 0  \\
0 & 0 & 0 & 1 & 0 & 0 & 0  \\
0 & 0 & 0 & 0 & 1 & 0 & 0  \\
0 & 0 & 0 & 0 & 0 & 1 & 0  \\
0 & 0 & 0 & 0 & 0 & 0 & 1  \\
0 & 1 & 0 & 0 & 0 & 0 & 0 
\end{pmatrix}
\end{equation}

\begin{figure}[hptb]
\begin{center}
\includegraphics[scale=0.3]{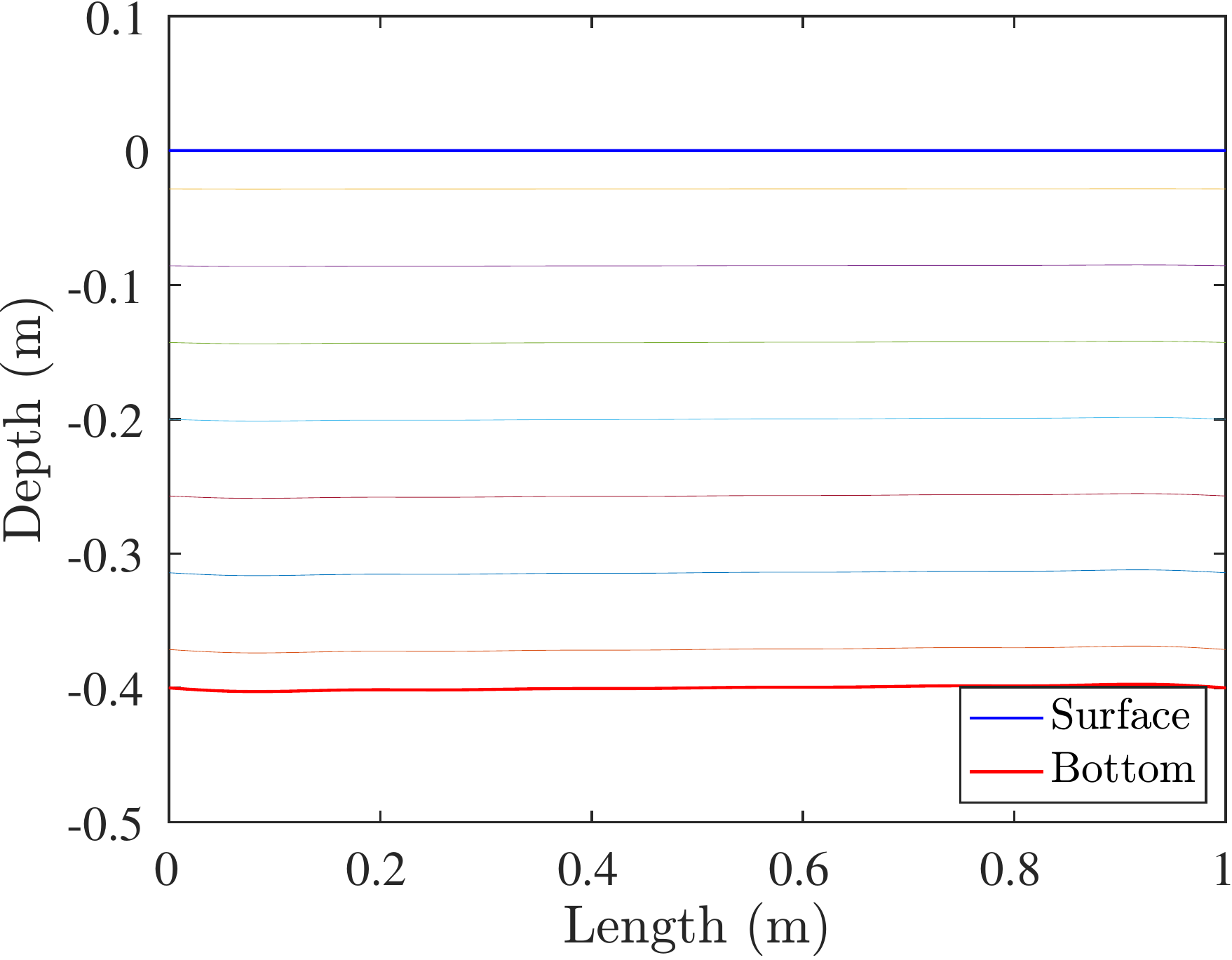}
\end{center}
\caption{The optimal topography and the associated trajectories for the permutation matrix~\eqref{eq:matL1}.}
\label{fig:top3}
\end{figure}

Finally we test for the length $L=\SI{1}{m}$. The optimal matrix $P_{\max}$ is given in~\eqref{eq:matL1} and the associated optimal topography is presented in Fig.~\ref{fig:top3}. The optimal value of $a$ is $[0.0017,0.0009,0.0006,0.0004,0.0003]$, and the two ratios defined in~\eqref{eq:r1}-\eqref{eq:r2} are $r_1=0.001\%$ and $r_2=3.453\%$. Note that the permutation strategy $P_{\max}$ is the same type as in~\cite[Equation 17]{Bernard02970756}.
\begin{equation}\label{eq:matL1}
P_{\max}=\begin{pmatrix}
1 & 0 & 0 & 0 & 0 & 0 & 0  \\
0 & 0 & 0 & 0 & 0 & 0 & 1  \\
0 & 0 & 0 & 0 & 0 & 1 & 0  \\
0 & 0 & 0 & 0 & 1 & 0 & 0  \\
0 & 0 & 0 & 1 & 0 & 0 & 0  \\
0 & 0 & 1 & 0 & 0 & 0 & 0  \\
0 & 1 & 0 & 0 & 0 & 0 & 0 
\end{pmatrix}
\end{equation}

\begin{figure}[hptb]
\begin{center}
\includegraphics[scale=0.3]{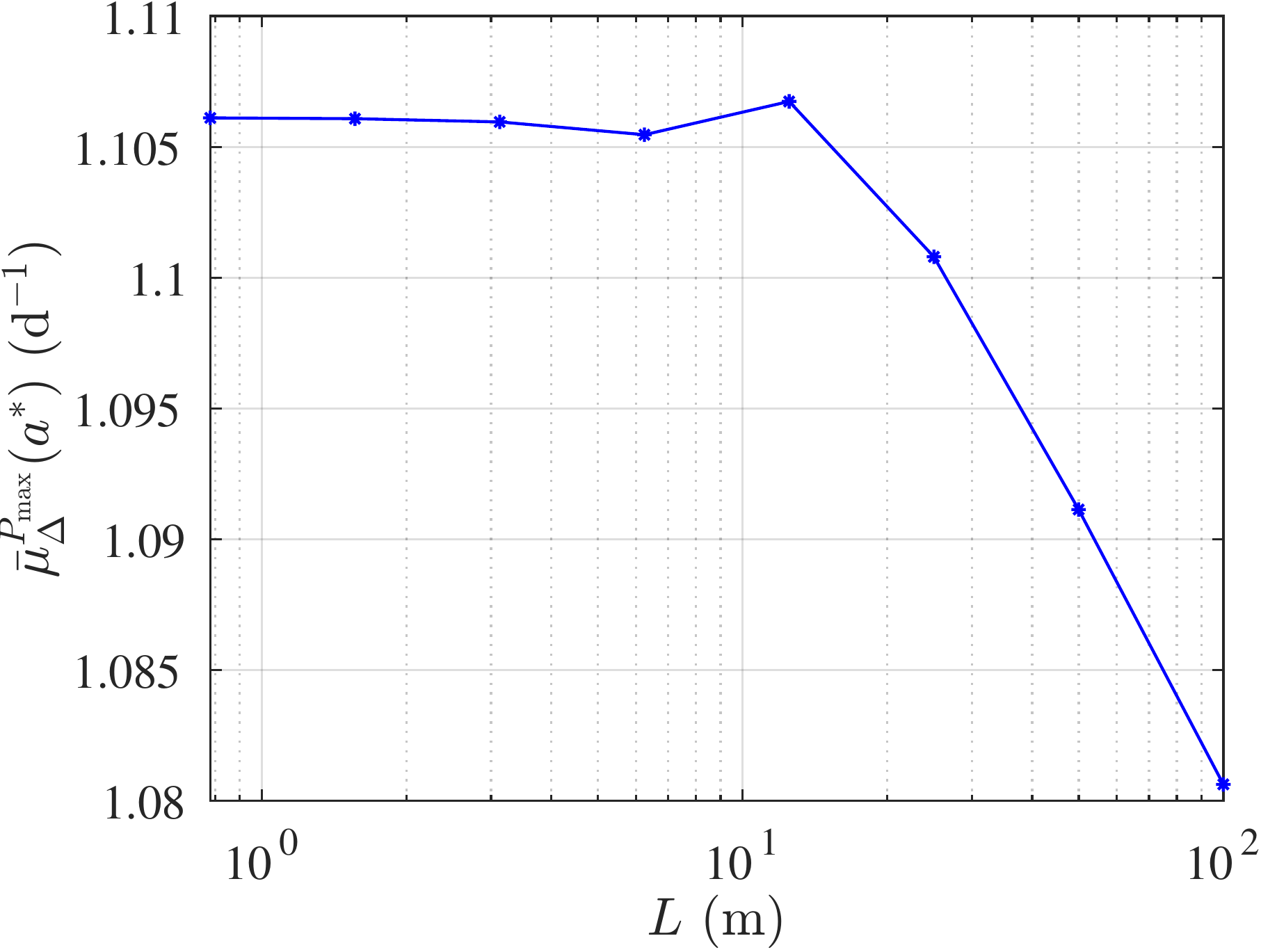}
\includegraphics[scale=0.3]{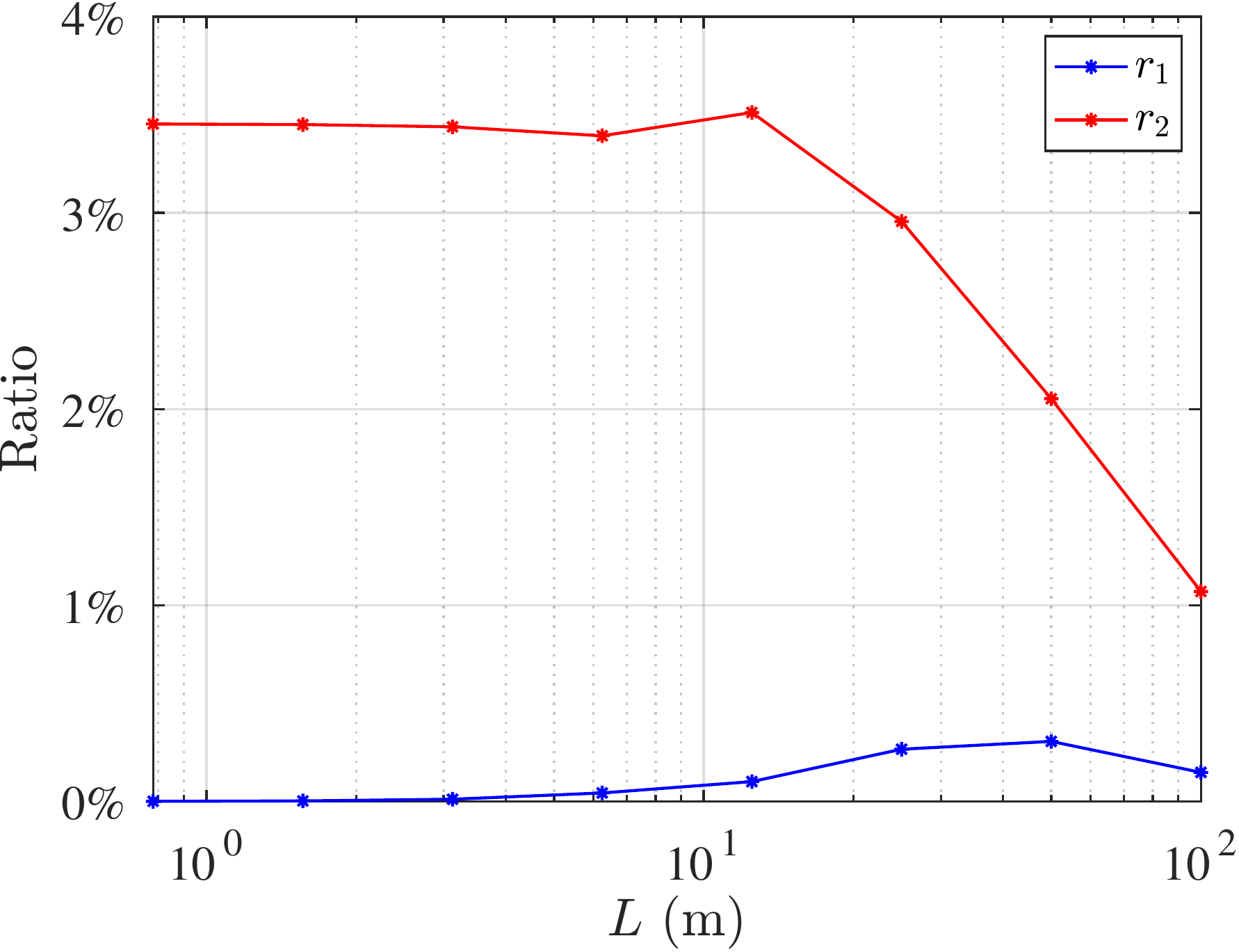}
\end{center}
\caption{The optimal value of the objective function $\bar \mu_{\Delta}$ (top) and the two ratios $r_1,r_2$ (bottom) for different values of length $L$.}
\label{fig:Lmur}
\end{figure}

As we observed in the three tests above, the length of raceway has a potential influence on the objective function and the gain, we then provide a test for different values of the length $L$. Fig.~\eqref{fig:Lmur} shows the objective function $\bar \mu_{\Delta}$ and the two ratios $r_1,r_2$ as a function of the length $L$. Note that the objective function decreases when $L$ increases except in the neighborhood of $L=\SI{12.5}{m}$, on the same time, we observe that the influence of topography is very limited comparing to the influence of the permutation strategies.

\begin{figure}[hptb]
\begin{center}
\includegraphics[scale=0.3]{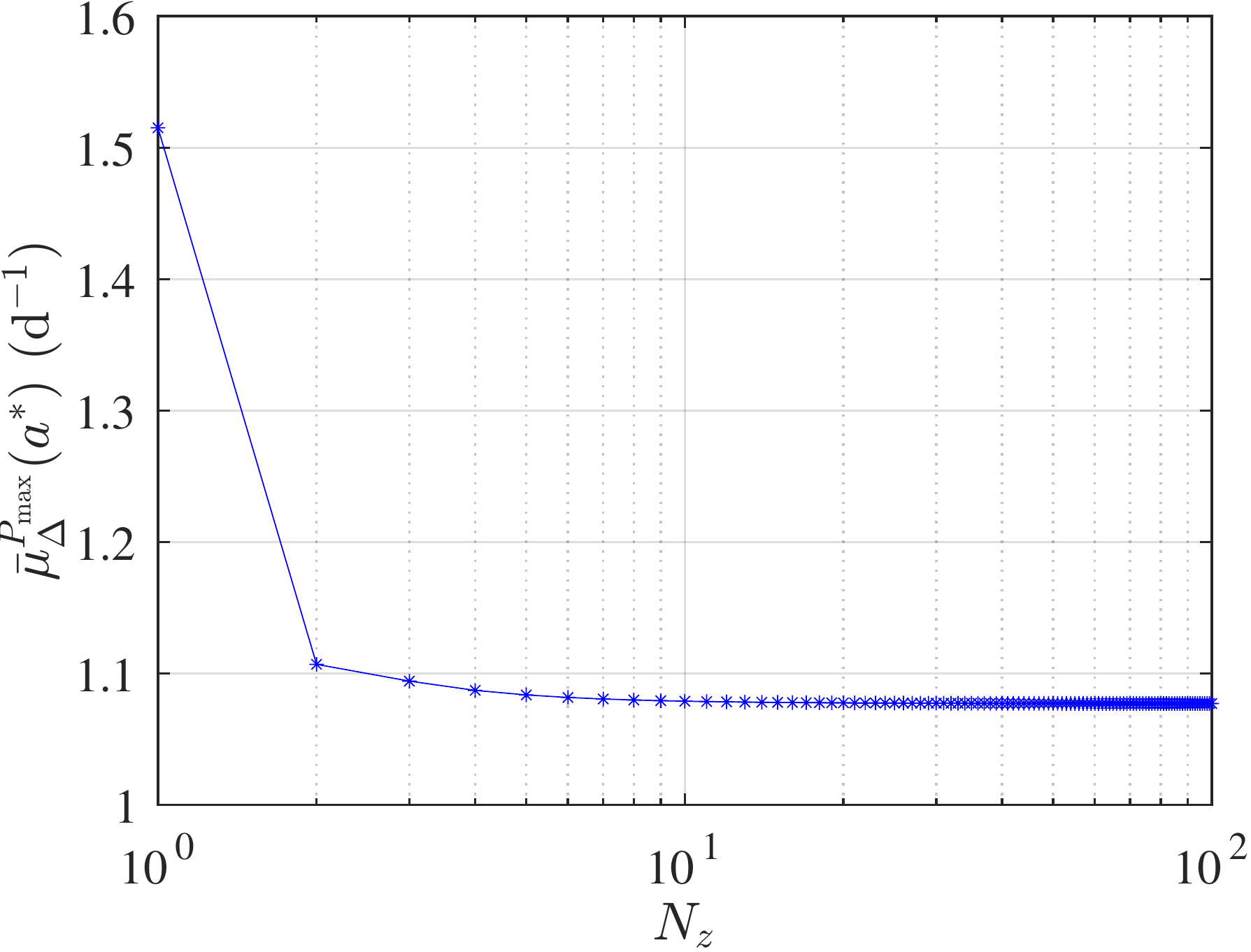}
\end{center}
\caption{The value of the objective function $\bar \mu_{\Delta}$ with a fix type of $P_{\max}$~\eqref{eq:matL100} for $N_z=[1,100]$.}
\label{fig:cv}
\end{figure}

As we mentioned above, we find the same type of the matrix as~\cite{Bernard02970756} which is for $N_z=11$ in the case $L=\SI{100}{m}$. Hence, we proceed further tests to investigate the convergence of the vertical discretization. Fig.~\eqref{fig:cv} presents the convergence of the objective function $\bar \mu_{\Delta}$ for the same type of the permutation strategy~\eqref{eq:matL100}.

\subsubsection{Varying volume}

\begin{figure}[hptb]
\begin{center}
\includegraphics[scale=0.3]{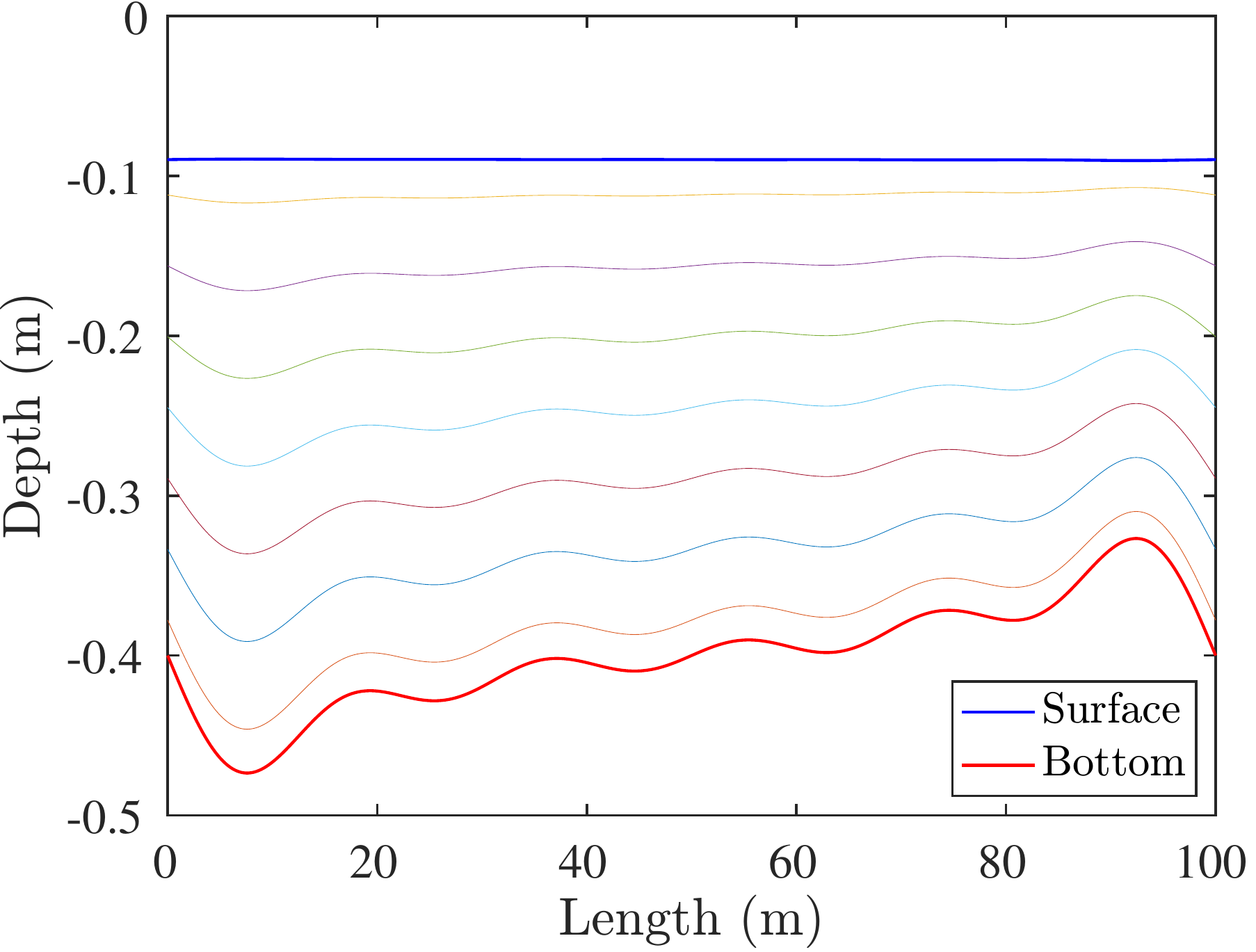}
\end{center}
\caption{The optimal topography and the associated trajectories for the permutation matrix~\eqref{eq:matVL100}.}
\label{fig:top4}
\end{figure}

We consider now that the volume can vary. Note that the volume related coefficient $a_0$ is also a parameter to be optimized. Let us define two ratios similar as~\eqref{eq:r1}-\eqref{eq:r2} to evaluate the efficiency of the permutation strategies,
\begin{align}
&\tilde r_1:=\frac{\Pi_{\Delta}^{P_{\max}}(\tilde a^*)-\Pi_{\Delta}^{P_{\max}}(\tilde a_f)}{\Pi_{\Delta}^{P_{\max}}(\tilde a_f)}, \label{eq:rt1}\\
&\tilde r_2:=\frac{\Pi_{\Delta}^{P_{\max}}(\tilde a^*)-\Pi_{\Delta}^{\I_{N_z}}(\tilde a_f)}{\Pi_{\Delta}^{\I_{N_z}}(\tilde a_f)},\label{eq:rt2}
\end{align}
where $\tilde a_f:=[\tilde a^*_0,0,\cdots,0]$ and $\tilde a^*_0$ is the optimal volume related value. Let us keep the truncated Fourier number and vertical discretization number settings as in previous test. The initial guess of $\tilde a$ is set to be a null vector with $\SI{40}{m^2}$ as initial value of the volume. For the length $L=\SI{100}{m}$, the optimal matrix $P_{\max}$ is given in~\eqref{eq:matVL100} and the associated optimal topography is presented in Fig.~\ref{fig:top4}. We obtain the optimal value $\tilde a^*=[0.3102,0.0328,0.0244,0.0173,0.0136,0.0129]$. The two ratios defined in~\eqref{eq:rt1}-\eqref{eq:rt2} are $\tilde r_1=0.686\%$ and $\tilde r_2=4.318\%$. 
\begin{equation}\label{eq:matVL100}
P_{\max}=\begin{pmatrix}
0 & 0 & 0 & 1 & 0 & 0 & 0  \\
0 & 0 & 0 & 0 & 0 & 1 & 0  \\
0 & 0 & 0 & 0 & 0 & 0 & 1  \\
0 & 0 & 0 & 0 & 1 & 0 & 0  \\
0 & 0 & 1 & 0 & 0 & 0 & 0  \\
0 & 1 & 0 & 0 & 0 & 0 & 0  \\
1 & 0 & 0 & 0 & 0 & 0 & 0 
\end{pmatrix}
\end{equation}

\begin{figure}[hptb]
\begin{center}
\includegraphics[scale=0.3]{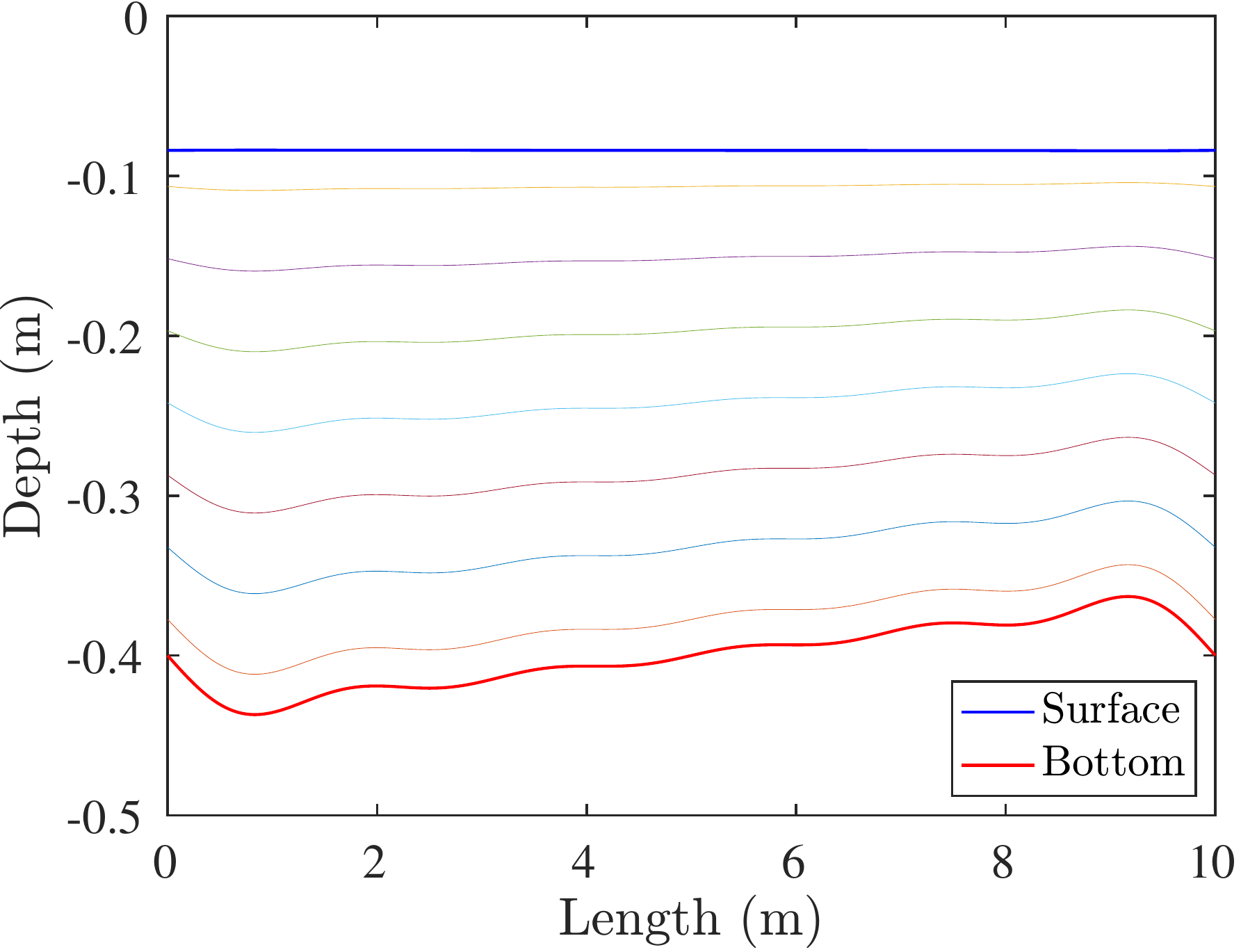}
\end{center}
\caption{The optimal topography and the associated trajectories for the permutation matrix~\eqref{eq:matVL10}.}
\label{fig:top5}
\end{figure}

For length $L=\SI{10}{m}$, the optimal matrix $P_{\max}$ is given in~\eqref{eq:matVL10} and the associated optimal topography is presented in Fig.~\ref{fig:top5}. The optimal value is $\tilde a^*=[0.3160, 0.0236, 0.0117, 0.0078, 0.0058, 0.0048]$, and the two ratios are $\tilde r_1=0.232\%$ and $\tilde r_2=12.299\%$. 
\begin{equation}\label{eq:matVL10}
P_{\max}=\begin{pmatrix}
0 & 0 & 0 & 0 & 0 & 0 & 1  \\
0 & 0 & 1 & 0 & 0 & 0 & 0  \\
0 & 0 & 0 & 1 & 0 & 0 & 0  \\
0 & 0 & 0 & 0 & 1 & 0 & 0  \\
0 & 0 & 0 & 0 & 0 & 1 & 0  \\
0 & 1 & 0 & 0 & 0 & 0 & 0  \\
1 & 0 & 0 & 0 & 0 & 0 & 0 
\end{pmatrix}
\end{equation}

\begin{figure}[hptb]
\begin{center}
\includegraphics[scale=0.3]{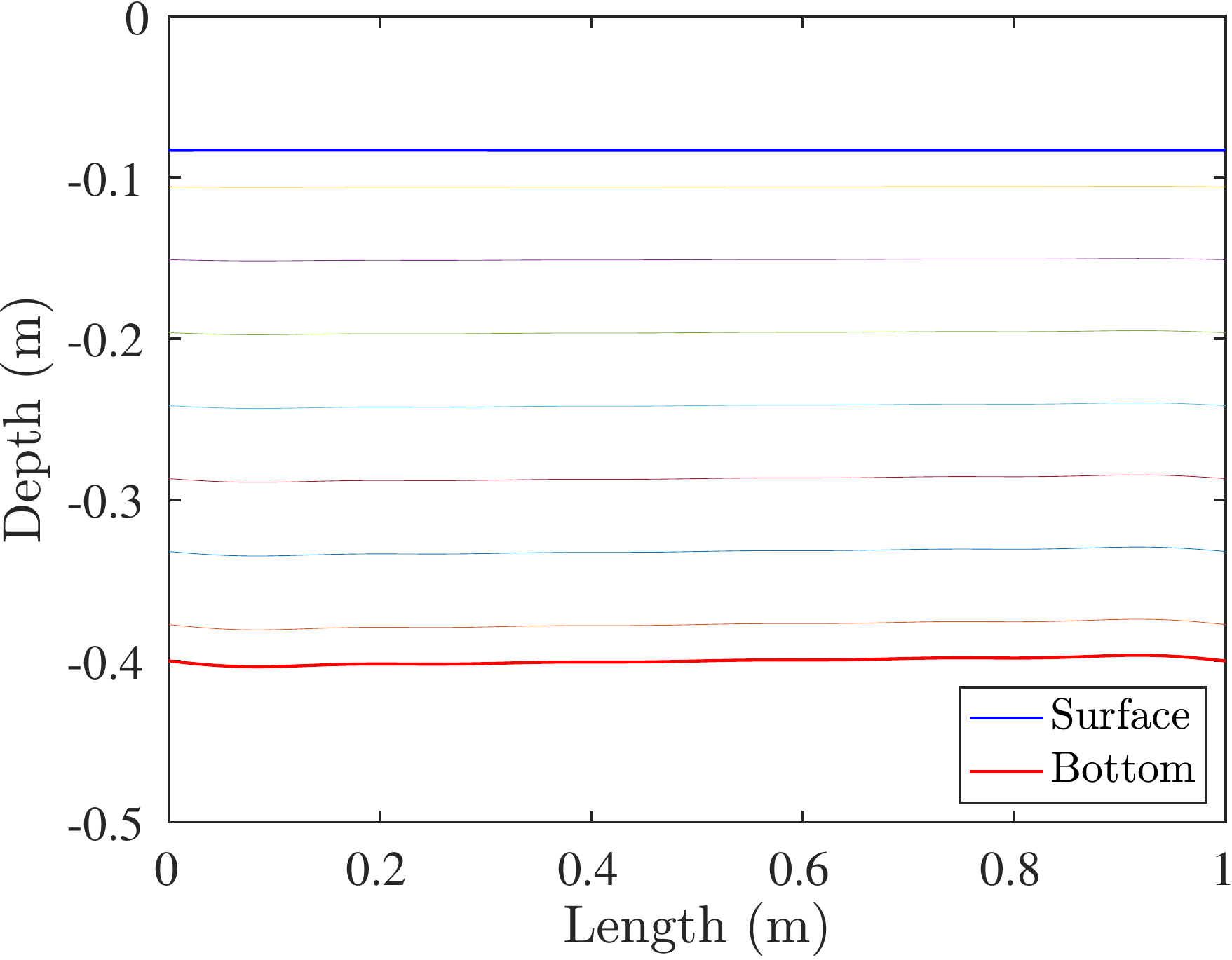}
\end{center}
\caption{The optimal topography and the associated trajectories for the permutation matrix~\eqref{eq:matVL1}.}
\label{fig:top6}
\end{figure}

For length $L=\SI{1}{m}$, the optimal matrix $P_{\max}$ is given in~\eqref{eq:matVL1} and the associated optimal topography is presented in Fig.~\ref{fig:top6}. The optimal value is $\tilde a^*=[0.3168, 0.0023, 0.0012, 0.0008, 0.0006, 0.0005]$, and the two ratios are $\tilde r_1=0.002\%$ and $\tilde r_2=12.714\%$. 
\begin{equation}\label{eq:matVL1}
P_{\max}=\begin{pmatrix}
0 & 0 & 0 & 0 & 0 & 0 & 1  \\
0 & 0 & 0 & 0 & 0 & 1 & 0  \\
0 & 0 & 0 & 0 & 1 & 0 & 0  \\
0 & 0 & 0 & 1 & 0 & 0 & 0  \\
0 & 0 & 1 & 0 & 0 & 0 & 0  \\
0 & 1 & 0 & 0 & 0 & 0 & 0  \\
1 & 0 & 0 & 0 & 0 & 0 & 0 
\end{pmatrix}
\end{equation}

\begin{figure}[hptb]
\begin{center}
\includegraphics[scale=0.3]{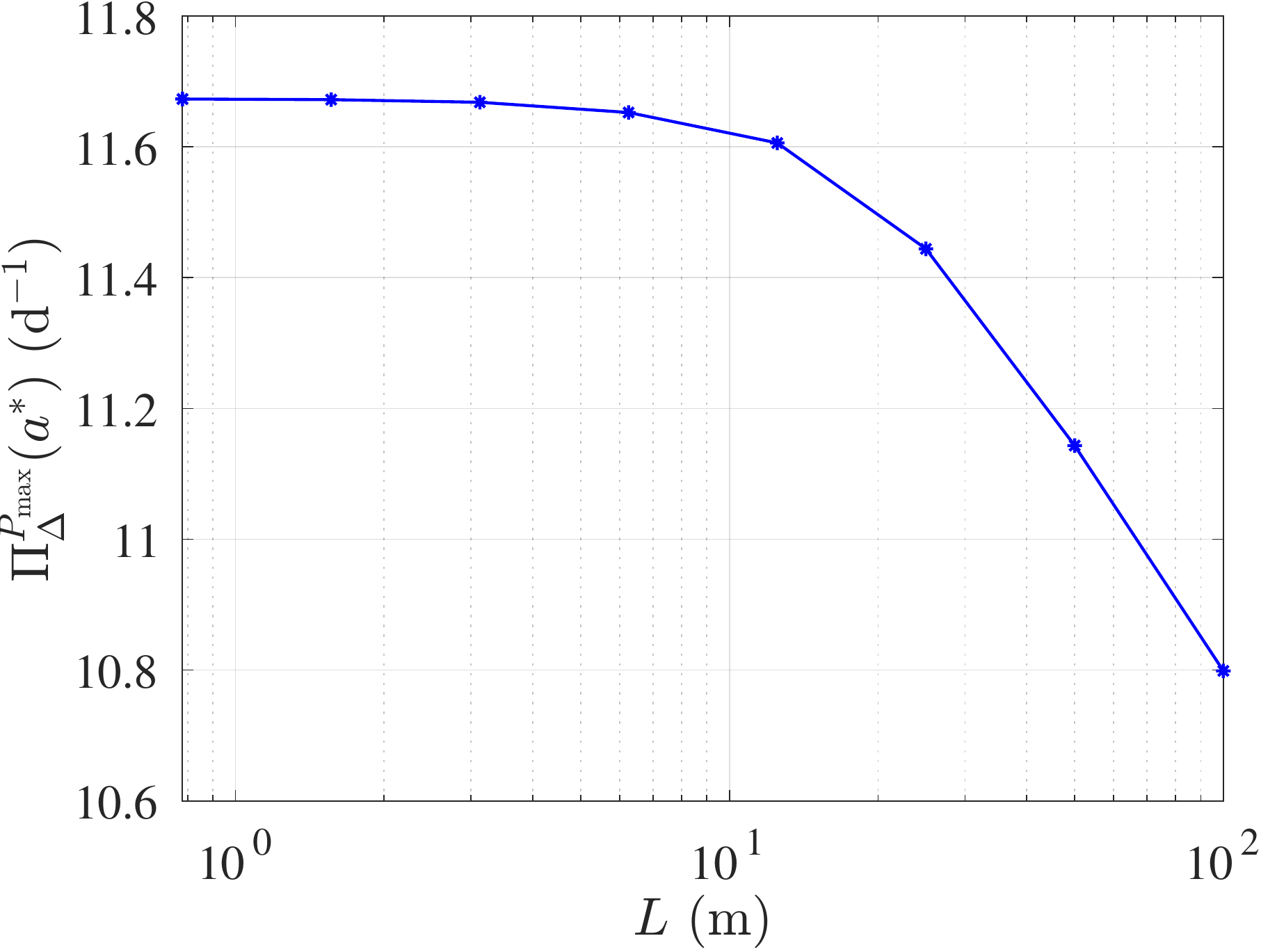}
\includegraphics[scale=0.3]{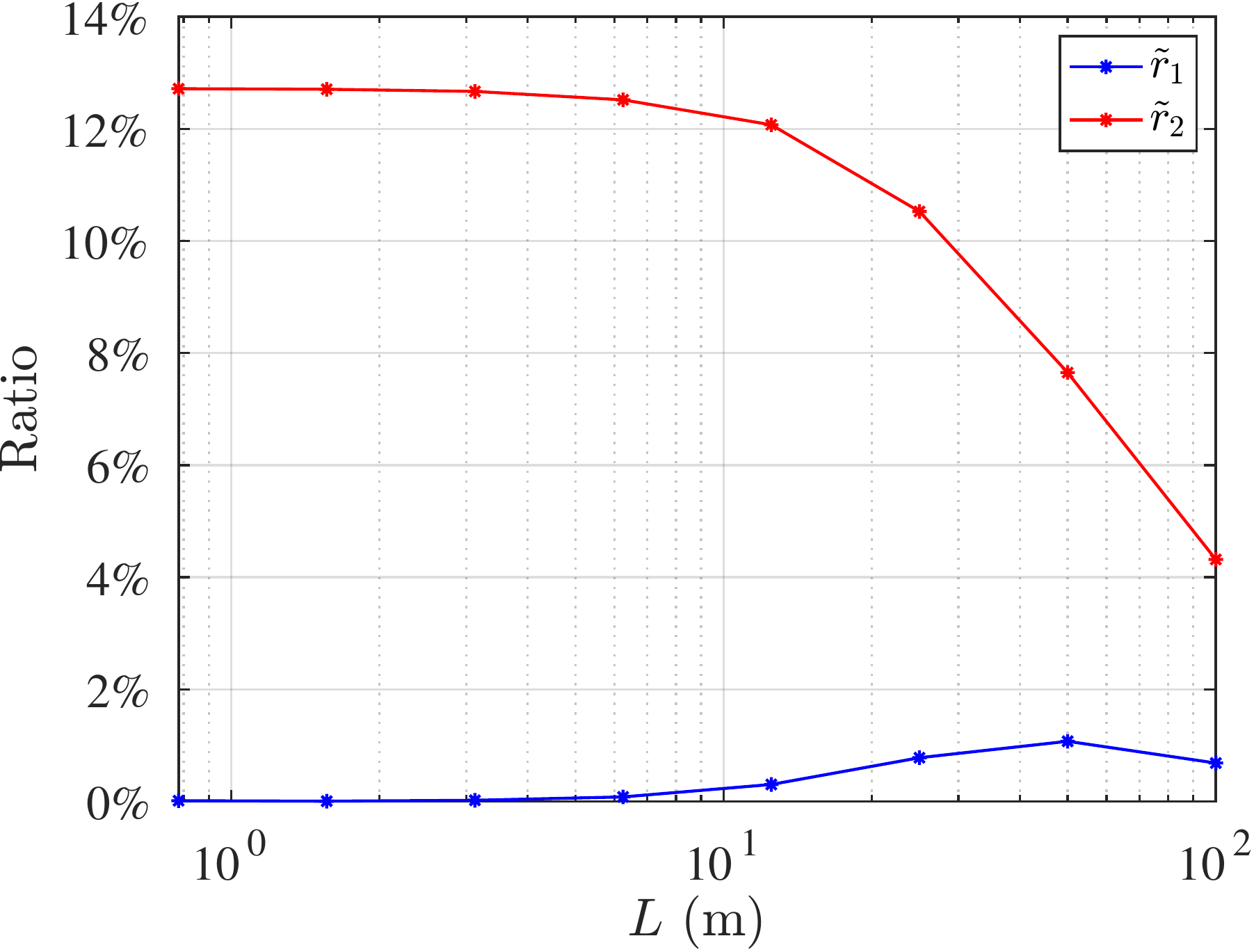}
\end{center}
\caption{The optimal value of the objective function $\Pi_{\Delta}$ (top) and the two ratios $\tilde r_1,\tilde r_2$ (bottom) for different value of length $L$.}
\label{fig:LPirtil}
\end{figure}

As shown experimentally in the previous three tests, the influence of the topographies remain limited, at the same time, non trivial permutation strategies $P_{\max}$ are obtained with different raceway length $L$, in particular these strategies are also different from the case with a fixed volume. Moreover, these strategies have a better improvement when the volume is also optimized. Fig.~\ref{fig:LPirtil} shows the objective function $\Pi_{\Delta}$ and the two ratios $\tilde r_1,\tilde r_2$ as a function of the length $L$. Note that the average growth rate $\Pi_{\Delta}$ increase monotonically when $L$ goes to 0. This \textit{flashing effect} corresponds to the fact that the algae exposed to high frequency flashing have a better growth. This phenomenon has already been reported in literature (e.g.~\cite{Bernard02970756,Lamare2018}).

\section{Conclusion}

Adapting the shape of the raceway to an original mixing system is an innovative strategy to boost the algal process productivity. To realize in practice the ideal mixing a system more elaborated than a paddle wheel is required. 

However, with the Han parameter considered for this species the gain stays limited and would not compensate a higher cost due to the more complicate design of the bottom topography and of the mixing device. It is possible that a higher gain is also obtained when leaving the laminar regime, but the energy dissipation in a turbulent regime would lead to a strong enhancement of the operating costs.

\bibliographystyle{plain}
\bibliography{ifacconf}

\end{document}